\documentclass{amsart}
\usepackage{graphicx}
\usepackage{amssymb}

\newtheorem{thm}{Theorem}
\newtheorem{stat}{Statement}

\newtheorem{prop}{Proposition}
\theoremstyle{definition}
\newtheorem{exam}{Example}

\theoremstyle{remark}
\newtheorem{rem}{Remark}
\numberwithin{equation}{section}

\newcommand{\set}[1]{\left\{#1\right\}}

\newcommand{\eps}{\varepsilon}

\title{On a class of C*-algebras determined uniquely by dual space.}
\author{Savchuk Yuri}%

\address{Kiev Taras Shevchenko National University, Faculty of
Cybernetics, Volody\-myr\-ska, 64, 01033, Kiev, Ukraine\\}

\email{ysavchuk@unicyb.kiev.ua}%

\begin{document}
\maketitle \vskip 0.5cm \noindent {\bf Mathematical Subject
Classification}: 55R10, 46L89, 46L99.
\begin{abstract}
Enveloping $C^*$-algebras for some finitely generated $*$-algebras
are considered. It is shown that all of the considered algebras
are identically defined by their dual spaces. The description in
terms of matrix-functions is given. \newline Keywords : algebraic
bundles, finitely presented $C^*$-algebras, dual spaces.
\end{abstract} \maketitle
\section*{Introduction}
One of the general approaches to description of $C^*$-algebras
consists in realization of them as "functions" on dual space. All
$C^*$-algebras having the irreducible representations in dimension
less or equal to some natural $n$, usually can be realized by
$n\times n$-matrix-functions with some boundary conditions. The
special case of such $C^*$-algebras is given by homogeneous
algebras - ones having all irreducible representations in the same
dimension $n$. The description of homogeneous $C^*$-algebras in
terms of fibre bundles was obtained in \cite{f},\cite{tt} (see
Statement 3 below). The analogous description for all
$C^*$-algebras having irreducible representations in different
dimensions less or equal to natural $n$ was presented in \cite{v}.

The homogeneous $C^*$-algebras whose spectra are tori $T_2,\ T_3$
(see \cite{dr}) and sphere $S^k$ (see \cite{ak}) were studied by
means of classification corresponding algebraic bundles. The
concrete example of non-homogeneous $C^*$-algebra is given by
$C^*$-algebra generated by the free pair of projections (see
\cite{vs}).

In the present paper we consider wide class of non-homogeneous
$F_{2n}$-$C^*$-algebras and give a realization of algebras from
this class in terms of matrix functions with boundary conditions.

\section{Preliminaries}
Let us recall some preparatory facts, definitions and notation,
which will be used below.

\subsection{Enveloping $C^*$-algebras}
In this paper we study finitely presented $C^*$-algebras. By
finitely presented $C^*$-algebra we mean an enveloping
$C^*$-algebra for finitely presented $*$-algebra.

{\bf Definition.} Let $\mathbf{A}$ be a $*$-algebra, having at
least one representation. Then a pair $(\mathcal{A},\rho)$ of a
$C^*$-algebra $\mathcal{A}$ and a homomorphism
$\rho:\mathbf{A}\longrightarrow\mathcal{A}$ is called an
enveloping pair for $\mathbf{A}$ if every irreducible
representation $\pi$:$\mathbf{A}\longrightarrow B(H)$ factors
uniquely through the $\mathcal{A}$, i.e. there exists precisely
one irreducible representation $\pi_1$ of algebra $\mathcal{A}$
satisfying $\pi_1\circ\rho=\pi$. The algebra $\mathcal{A}$ is
called an enveloping for $\mathbf{A}$.

\subsection{$F_{2n}$ - algebras}

Since all $C^*$-algebras considered below are examples of
$F_{2n}$-$C^*$-algebras, we recall here the definition and basic
properties of algebras satisfying the $F_{2n}$ polynomial
identity. Let $F_n$ denotes the following polynomial of degree $n$
in $n$ non-commuting variables:
$$F_n(x_1,x_2,\ldots,x_n)=\sum_{\sigma\in
S_n}(-1)^{p(\sigma)}x_{\sigma(1)}\ldots x_{\sigma(n)},
$$
where $S_n$ is the symmetric group of degree $n$, $p(\sigma)$ is
the parity of a permutation $\sigma\in S_n$. We say that an
algebra $A$ is an algebra with $F_n$ identity if for all
$x_1,\ldots,\ x_n\in A$, we have $F_n(x_1,\ldots x_n)=0$. The
Amitsur-Levitsky theorem says that the matrix algebra
$M_n(\mathbb{C})$ is an algebra with $F_{2n}$ identity. Moreover,
$C^*$-algebra $A$ has irreducible representations of dimension
less or equal to $n$ iff $A$ satisfies the $F_{2n}$ condition (see
\cite{os}). An important class of $F_{2n}$-$C^*$-algebras is
formed by $n$-homogeneous algebras - ones having all irreducible
representations of dimension $n$. The simplest example of
$n$-homogeneous $C^*$-algebra is $2\times 2$-matrix algebra over
$C(X)$, where $X$ is a compact Hausdorff space - it is so-called
trivial algebra. In this case the dual space (space of primitive
ideals, see \cite{d}) is canonically homeomorphic to $X$. Simple
example of non-homogeneous $C^*$-algebra is given by "fixing" the
homogeneous algebra above in a point $x_0\in X$. Namely, consider
$$A(x_0)=\set{f\in C(X\longrightarrow M_2(\mathbb{C}))|
f(x_0)\in M_1(\mathbb{C})\oplus M_1(\mathbb{C})}.$$

The dual space $P(A(x_0))$ of algebra $A(x_0)$ is $T_2\cup T_1$,
where $T_2$ is the space of irreducible $2$-dimensional
representations, homeomorphic to $X\backslash\set{x_0}$ and
$T_1=\set{x_0',x_0''}$ is the space of $1$-dimensional
representations. Here the bases of neighborhoods of the points
$x_0'$ and $x_0''$ are the same as for the $x_0$ in $X$, i.e. we
have the space satisfying $T_0$ but not $T_1$ separability axiom
(such spaces sometimes are called quasicompact).

The description above can be easily generalized, for example
taking $n\geq 2$:
$$A(x_0)=\set{f\in C(X\longrightarrow M_n(\mathbb{C}))|
f(x_0)\in M_{n_1}(\mathbb{C})\oplus\ldots\oplus
M_{n_k}(\mathbb{C})},$$ or "fixing" a homogeneous $C^*$-algebra in
a few points. Studying some examples below, we will describe dual
spaces in terms of algebras
\begin{gather*}
A_n(m_1^{(1)},\ldots,m_{k_1}^{(1)};m_1^{(2)},\ldots,m_{k_2}^{(2)};m_1^{(3)},\ldots,m_{k_3}^{(3)})=\\
=\{f\in C(S^2\longrightarrow M_n(\mathbb{C}))| f(x_i)\in
M_{m_1^{(i)}}(\mathbb{C})\oplus\ldots\oplus
M_{m_{k_i}^{(i)}}(\mathbb{C}),\ x_i\in S^2,\ i=1,2,3\},
\end{gather*}
where $x_1,x_2,x_3\in S^2$ are fixed points (it can be checked
that the different triples $x_1,x_2,x_3$ give isomorphic
algebras).

Finitely presented $C^*$-algebras are mostly non-homogeneous.
Below we present some important examples of such algebras.

\begin{exam}
Main examples of the paper concern $*$-algebras ($C^*$-algebras)
associated with extended Dynkin diagrams
$\widetilde{D_4},\widetilde{E_6},\widetilde{E_7},\widetilde{E_8}$,
see \cite{vms}, for details. Below we recall how these algebras
are defined by generators and basic relations:
\newline
$C^*(\widetilde{D_4})= C^*\langle
p_1,p_2,p_3,p_4|p_1+p_2+p_3+p_4=2e,\
p_i^*=p_i=p_i^2,i=\overline{1,4}\rangle, $
\newline
$ C^*(\widetilde{E_6})=
C^*\langle p_1,p_2,q_1,q_2,r_1,r_2|p_1+2p_2+q_1+2q_2+r_1+2r_2=3e,\\
p_i^*=p_i=p_i^2,q_i^*=q_i=q_i^2,r_i^*=r_i=r_i^2,i=\overline{1,2},
p_jp_k=q_jq_k=r_jr_k=0,j\neq k\rangle,$
\newline
$
C^*(\widetilde{E_7})=C^*\langle p_1,p_2,p_3,q_1,q_2,q_3,r|p_1+2p_2+3p_3+q_1+2q_2+3q_3+2r=4e,\\
p_i^*=p_i=p_i^2,q_i^*=q_i=q_i^2,i=\overline{1,3},r^*=r=r^2,p_jp_k=q_jq_k=0,j\neq
k\rangle, $
\newline
$
C^*(\widetilde{E_8})=C^*\langle p_1,p_2,q_1,q_2,q_3,q_4,q_5,r|2p_1+4p_2+q_1+2q_2+3q_3+4q_4+5q_5+3r=6e,\\
p_i^*=p_i=p_i^2,q_j^*=q_j=q_j^2,r^*=r=r^2,i=\overline{1,2},j=\overline{1,5},p_lp_m=q_lq_m=0,l\neq
m\rangle.$

Irreducible representations of the algebras above were classified
in the series of works (see \cite{os}, \cite{m}, \cite{o},
\cite{m1}) and the description of dual spaces follows easily from
the classification (see for example \cite{yu}). The dual spaces of
$C^*$-algebras presented above are the following:

\begin{itemize}
\item for $C^*(\widetilde{D_4})$ the dual space is the same as for
$A_2(1,1;\ 1,1;\ 1,1)$

\item for $C^*(\widetilde{E_6})$ the dual space is the same as for
$A_3(2,1;\ 1,1,1;\ 1,1,1)$

\item for $C^*(\widetilde{E_7})$ the dual space is the same as for
$A_4(2,2;\ 2,1,1;\ 1,1,1,1)$

\item for $C^*(\widetilde{E_8})$ the dual space is the same as for
$A_6(3,3;\ 2,2,1,1;\ 1,1,1,1,1,1)$.

\end{itemize}

\end{exam}

\begin{exam}
Generalizing the relations defined by $\widetilde{D_4}$ one can
consider the $*$-algebra:
\begin{gather*}
P_{\alpha,\beta}=\mathbb{C}\langle p_1,p_2,p_3,p_4|
\alpha(p_1+p_2)+\beta(p_3+p_4)=I\rangle,
\end{gather*}
where $\alpha,\beta>0,\alpha+\beta=1$. The case $P_{\frac 1 2,
\frac 1 2}$ corresponds to $\widetilde{D_4}$. When
$\alpha\neq\beta$, the dual space for $C^*(P_{\alpha,\beta})$ is
the same as for algebra:
\begin{gather*}
\{f\in C(S^2\longrightarrow M_2(\mathbb{C}))|f(x_1),\ f(x_2)\in
M_1(\mathbb{C})\oplus M_1(\mathbb{C}),\ x_1,\ x_2\in S^2\},
\end{gather*}
it can be checked that the different pairs $(x_1,x_2)$ give
isomorphic algebras, see \cite{kir} for more details.
\end{exam}

\begin{exam}\label{com}
We can consider much more complicated relations using the algebras
$P_{\alpha,\beta}$ from previous example. Construct the next
algebra:
\begin{gather*}
\overline{P}_{\eps}=C^*\langle
p_1,p_2,p_3,p_4,\alpha,\beta| \alpha(p_1+p_2)+\beta(p_3+p_4)=I,\alpha=\alpha^*,\beta=\beta^*\\
[\alpha,p_i]=[\beta,p_i]=[\alpha,\beta]=0,\ i=\overline{1,4},\
\alpha,\beta\geq\eps I,\ \alpha+\beta=I\rangle,
\end{gather*}
where $0<\eps<\frac 1 2$. Using the previous example one can
describe the dual space for the $\overline{P}_\eps$. We will see
later that all algebras $\overline{P}_\eps,\ \eps\in (0,1/2)$ are
isomorphic.
\end{exam}

\begin{exam}\label{free}
The group algebra for $G=\mathbb{Z}_2\ast \mathbb{Z}_2$ gives an
example of $F_4$-algebra corresponding to infinite discrete group.
$$C^*(\mathbb{Z}_2\ast \mathbb{Z}_2)=C^*\langle p_1,p_2|p_k^2=p_k=p_k^*,\
k=1,2\rangle:=\mathcal{P}_2.
$$
\end{exam}

\subsection{Fibre bundles}
All the necessary information about fibre bundles reader can find
in \cite{mi}.

From the general theory of fibre bundles it is known that
locally-trivial bundle is defined by its 7-tuple\newline
$(E,B,F,G,p,\set{V_i}_{i\in I},\set{\phi_i}_{i\in I})$, where $E$
is a bundle space, $B$ is a base, $F$ is a fibre, $G$ is a
structure group, $p:E\to B$ is a projection, $V_i$ is an open
covering of base $B$ and $\phi_{i}:V_i\times F\to p^{-1}(V_i)$ are
coordinate homeomorphisms. Note that in Vasil'ev's paper \cite{v}
coordinate and skew products are considered instead of locally
trivial bundles. Coordinate product is a 7-tuple
$(E,B,F,G,p,\set{V_i}_{i\in I},\set{\phi_{ji}}_{i,j\in I})$, where
$\phi_{ji}=\phi_j^{-1}\circ\phi_i|_{V_i\cap V_j}:(V_i\cap
V_j)\times F\rightarrow (V_i\cap V_j)\times F$ are sewing maps.
Every $\phi_{ji}$ is uniquely determined by continuous map (which
we also call sewing map) $\overline{\phi}_{ji}:V_i\cap V_j\to G$
by formulae: $\phi_{ji}(x,f)=(x,\overline{\phi}_{ji}(x)\circ f),\
\forall x\in V_i\cap V_j$.

In this paper we consider only algebraic bundles, i.e. bundles
with $F\simeq M_n(\mathbb{C})$ and $G\subseteq\underline{U(n)}$,
where $\underline{U(n)}$ is a group of inner automorphisms of
$M_n(\mathbb{C})$, i.e. group of unitary matrices $U(n)$
factorized by it's center. Group $\underline{U(n)}$ is supposed to
act on the $M_n(\mathbb{C})$ in the following way. Let
$\hat{g}\in\underline{U(n)}$ and $g\in U(n)$ is corresponding
matrix, then $\hat{g}(m)=g^*m g,\ m\in M_n(\mathbb{C})$. If $G$ is
a subgroup of $U(n)$, we denote by $\underline{G}$ the quotient of
$G$ by it's center. Sometimes we suppose sewing maps
$\overline{\phi}_{ji}$ to have values in $U(n)$, instead of
$\underline{U(n)}$, one has to take a composition of
$\overline{\phi}_{ji}$ with the canonical map of $U(n)$ onto
$\underline{U(n)}$.

The following statements on triviality will be used below.

\begin{stat}$($see \cite{hu}$)$
Every locally-trivial bundle over $n$-dimensional disk $D^n$ is
isomorphic to trivial.
\end{stat}

\begin{stat}$($see \cite{ak}$)$
Every algebraic bundle over some compact subset of $\mathbb{C}$ is
isomorphic to trivial.
\end{stat}

In this article by $\mathcal{B}(X,F,G)$ we denote some bundle with
base $X$, fibre $F$ and group $G$. The algebra of sections of
bundle $\mathcal{B}$ is denoted by $\Gamma(\mathcal{B})$.

\subsection{Structure of $F_{2n}$-$C^*$-algebras}

Next statement (Fell-Tomiyama-Takesaki theorem) gives a
description of $n$-homogeneous $C^*$-algebras in terms of
algebraic bundles.

\begin{stat}$($\cite{tt},\cite{v}$)$
For every $n$-homogeneous $C^*$-algebra $A$ there exists a bundle
$\mathcal{B}(P(A),M_n(\mathbb{C}),\underline{U(n)})$, such that
$A\simeq \Gamma(\mathcal{B})$, where $P(A)$ is space of all pairly
non-equivalent irreducible representations.
\end{stat}

Analogous result for non-homogeneous $F_{2n}$-$C^*$-algebras was
obtained in Vasil'ev's work \cite{v}. The following statement is
the straightforward corollary from the main result of \cite{v}.

\begin{stat}\label{vas}
Let $A$ be $F_{2N}$-$C^*$-algebra, having finite number of
irreducible representations in dimensions $1,2,\ldots,N-1$. Then
there exist:
\begin{enumerate}
\item finite sets $X_1,\ldots,X_{N-1}$, compact space
$\overline{X}$ and open dense in $\overline{X}$ subspace $X$, such
that $\overline{X}\backslash X$ is finite,

\item a formal decomposition for all $x\in\overline{X}\backslash
X$ :
$$x=\sum_{k=1}^r d_k x_k,$$
where $d_k$ are natural numbers, $x_k\in X_{n_k}$ and $\sum d_k
n_k=N$, such that two points from $\overline{X}\backslash X$
coincide iff they have the same decomposition,

\item bundle
$\mathcal{B}=\mathcal{B}(X,M_N(\mathbb{C}),\underline{U(N)})$
having $\underline{n_1U(d_1)\times\ldots\times n_rU(d_r)}$ as the
structure group at the point $x\in \overline{X}\backslash X,\
x=d_1x_1+\ldots+d_rx_r$ (i.e. there exists a neighborhood
$O_x\subseteq \overline{X}$, such that the structure group of
bundle $\mathcal{B}$ reduces to
$\underline{n_1U(d_1)\times\ldots\times n_rU(d_r)}$ at the
subspace $O_x\cap X$),
\end{enumerate}
such that algebra $A$ can be realized as an algebra of $N$-tuples
$a=(a^{(1)},\ldots, a^{(N)})$, where $a^{(i)}:X_i\longrightarrow
M_i(\mathbb{C}),\ i=\overline{1,N-1}$ are some functions and
$a^{(N)}\in\Gamma(\mathcal{B})$ is a section such that:

\begin{gather*}\lim_{x\to\sum_{k=1}^r d_k x_k}a^{(N)}(x)=
\begin{pmatrix}
  \mathbf{1}_{d_1}\otimes a^{(n_1)}(x_1) & & & \mbox{{\huge $0$}}\\
    & \ddots \\
  \mbox{{\huge $0$}} & & & \mathbf{1}_{d_r}\otimes a^{(n_r)}(x_r)
\end{pmatrix}.
\end{gather*}

\end{stat}

\subsection{Bundles over 2-sphere}
Recall, that fibre bundles over the sphere $S^k$ are naturally
classified by elements of the homotopic group $\pi_{k-1}(G)$,
where $G$ is structure group of bundle (see for example \cite{mi}
or \cite{hu}).

Bundles over 2-dimensional sphere with fibre $M_n(\mathbb{C})$ are
classified in \cite{ak}, in this case $\pi(\underline{U(n)})\simeq
\mathbb{Z}_n$. Let us give a sketch of this classification.

Consider an atlas on 2-sphere, containing 2 charts : upper and
lower closed half spheres ($S^2_+=W_1$ and $S^2_-=W_2$
correspondingly). As sewing map $\overline{\phi}_{1,2}$ take
continuous map $V':S^1\longrightarrow \underline{U(n)}$. Construct
a bundle $\mathcal{B}'$ with atlas $\set{W_1,\ W_2}$ and sewing
map $V'$. Two bundles $\mathcal{B}'$ and $\mathcal{B}''$ defined
by sewing maps $V'$ and $V''$ are isomorphic iff $ind\ V'-ind\
V''=n l,\ l\in \mathbb{Z}$, where $ind\ V=(2\pi)^{-1}[\arg \det
V(t)]_{S^1}$ is the winding number of the function $V$ with
respect to the circle. Define sewing maps:
\begin{gather}\label{Vk}
V_k:S^1\longrightarrow U(n), z\mapsto diag(z^k,1,\ldots,1),\
k=\overline{{0,n-1}}.
\end{gather}

Then the bundles $\mathcal{B}_{n,k},\ k=\overline{0,n-1}$, defined
by $V_k$ are canonical representatives of all isomorphism classes.

\section{Enveloping $C^*$-algebras, an examples.}

To give a description of $C^*$-algebras for examples from Section
1, we need the following

\begin{thm}
Let
\begin{gather*}D^l(m_1,\ldots,m_k)=\{f\in
C(D^l\longrightarrow M_{m_1+\ldots+m_k}(\mathbb{C}))|\\
f(0,\ldots,0)\in M_{m_1}(\mathbb{C})\oplus\ldots\oplus
M_{m_k}(\mathbb{C})\},
\end{gather*}

Every $C^*$-algebra $A$ having the dual space same as
$D^l(m_1,\ldots,m_k)$ is isomorphic to $D^l(m_1,\ldots,m_k)$.
\end{thm}
\begin{proof}
We give proof when $l=2,\ k=2,\ m_1=m_2=1$, i.e.
$D^l(m_1,\ldots,m_k)=D^2(1,1)$ and realize $D^2$ as the unit disk
in complex plane. For the general case proof is analogous. Let us
introduce some notaion:
$$
D_\pm^2(1,1)=\set{f\in C(D^2_\pm\longrightarrow
M_2(\mathbb{C}))|f(0,0)\in M_1(\mathbb{C})\oplus M_1(\mathbb{C})},
$$
where $D^2_+=\set{z\in D^2| \Im m(z)\geq 0}$ and $D^2_-=\set{z\in
D^2|\Im m(z)\leq 0}$, and let $A_+$ be $C^*$-algebra having dual
space same as for $D^2_+$. We have to prove that $A_+\simeq
D_+^2$. The statement \ref{vas} implies that there exists a bundle
$\mathcal{B}_+=\mathcal{B}_+(D^2_+\backslash\{0\},M_2(\mathbb{C}),\underline{U(2)})$,
having $\underline{U(1)\times U(1)}$ as the structure group at the
point $0$, such that $A_+$ is isomorphic to the algebra of pairs:

$$\set{(a',a'')|a'=a'(x)\in\Gamma(\mathcal{B_+}),a''\in M_1(\mathbb{C})
\oplus M_1(\mathbb{C}),\lim_{x\to 0}a'(x)=a''}.$$

We suppose $\mathcal{B_+}$ to be assigned by charts
$\set{V_k}_{k\in \mathbb{N}}$, where
$$V_k=\set{z\in D^2_+,\ \frac 1 {2^k}\leq |z|\leq \frac 1 {2^{k+1}}}.$$

Sewing maps $\overline{\phi}_{k,k+1}$ of the bundle
$\mathcal{B}_+$ are $\set{\mu_k}_{k\in \mathbb{N}}$, where
$\mu_k:V_k\cap V_{k+1}\rightarrow \underline{U(2)}$ are some
continuous functions.

The condition for $\mathcal{B}$ to have $\underline{U(1)\times
U(1)}$ as the structure group at the point $0$ implies that there
exists $p\in \mathbb{N}$ such that
$\mu_r(z)\in\underline{U(1)\times U(1)},\forall r>p$. The fact
that $\underline{U(2)}$ and $\underline{U(1)\times U(1)}$ are
connected groups allows us to construct continuous maps
$$\overline{\mu_k}:V_k\rightarrow \underline{U(2)},\ k\leq p$$ and
$$\overline{\mu_k}:V_k\rightarrow \underline{U(1)\times U(1)},k>p$$
such that
\begin{gather*}
\overline{\mu_k}(z)=\left\{%
\begin{array}{ll}
    e, & z\in V_{k-1}\cap V_k; \\
    \mu_k(z), & z\in V_k\cap V_{k+1}. \\
\end{array}%
\right.
\end{gather*}

Let us write sections of $\mathcal{B}_+$ explicitly:
$$\Gamma(\mathcal{B}_+)\ni a'\longleftrightarrow\set{a'_k}_{k\in \mathbb{N}},
$$
where $a'_k:V_k\rightarrow M_2(\mathbb{C})$ are continuous
functions satisfying the condition of compatibility
$a'_{k+1}(z)=\mu_k(z)\circ(a'_k(z)),\ z\in V_k\cap V_{k+1}$. So we
write elements $(a',a'')\in A_+$ as $(\{a'_k\}_{k\in
\mathbb{N}},a'')$.

Let us consider the map:
$$A_+\ni(\{a'_k(z)\}_{k\in\mathbb{N}},a'')\mapsto
(\{\overline{\mu}_k(z)\circ a'_k(z)\}_{k\in\mathbb{N}},a'').
$$
One can check that the map is the isomorphism of algebras $A_+$
and $D^2_+(1,1)$. Analogously the isomorphism $A_-\simeq
D_-^2(1,1)$ can be proved.

Let us consider bundle
$\mathcal{B}(D^2\backslash\{0\},M_2(\mathbb{C}),\underline{U(2)})$
having $\underline{U(1)\times U(1)}$ as the structure group at the
point $0$, such that:
$$A\simeq\set{(a',a'')|a'=a'(x)\in\Gamma(\mathcal{B}),a''\in M_1(\mathbb{C})
\oplus M_1(\mathbb{C}),\lim_{x\to 0}a'(x)=a''}.$$

The fact that $A_+$ and $A_-$ are "trivial" allows us to suppose
$\mathcal{B}$ to have only two charts: $W_1=D^2_+\backslash\{0\}$
and $W_2=D^2_-\backslash\{0\}$ and one sewing map
$\mu_{1,2}:[-1,1]\backslash \{0\}\rightarrow \underline{U(2)}$
such that $\exists\eps_0>0:\forall z\in (-\eps_0,0)\cup(0,\eps_0)\
\mu(z)\in\underline{U(1)\times U(1)}$. As before, by connectivity
of $\underline{U(2)}$ and $\underline{U(1)\times U(1)}$ we
construct continuous map $\overline{\mu}:W_1\rightarrow
\underline{U(2)}$ possessing a values from $\underline{U(1)\times
U(1)}$ in $O(\eps_0,0)=\{z\in D^2_+,0<|z|<\eps_0\}$ such that
$\forall z\in W_1\cap W_2 : \overline{\mu}(z)=\mu_{1,2}(z)$.

As before we have an isomorphism given by $\overline{\mu}$:
$$A\ni(\{a'_1(z),a'_2(z)\},a'')\mapsto (\{\overline{\mu}(z)a'_1(z),a'_2(z)\},a'')\in D^2(1,1).
$$
The proof is completed.
\end{proof}

\begin{rem}
One can regard the theorem above as a certain generalization of
the Statement 1. Indeed we have a $C^*$-algebra identically
defined by it's dual space.

\end{rem}
As before, let
$$\sum_{j=1}^k m_j=n,\ m_j\in\mathbb{N},j=\overline{1,k}$$ and set
$$B_{n,p}(m_1,\ldots,m_k)=\{f\in \Gamma(\mathcal{B}_{n,p}) |
f(x)\in M_{m_1}(\mathbb{C})\oplus\ldots\oplus
M_{m_k}(\mathbb{C})\},$$ where $x\in S^2$ is fixed point (we
consider $x$ to be covered only by one chart).

Next theorem is basic to study the structure of $C^*$-algebras
corresponding to Dynkin diagrams.

\begin{thm}
If $l_1\equiv l_2\ (mod (m_1,\ldots,m_k)\ )$, then
\begin{gather}\label{bnl12}B_{n,l_1}(m_1,\ldots,m_k)\simeq
B_{n,l_2}(m_1,\ldots,m_k)
\end{gather}

(here $(m_1,\ldots,m_k)$ denotes the greatest common divider of
$m_1,\ldots,m_k$)
\end{thm}

\begin{proof}
Let us realize elements of algebras $B_{n,p}(m_1,\ldots,m_k)$ as
sections of corresponding bundles $\mathcal{B}_{n,p}$ (here
$p=l_1$ or $p=l_2$):
\begin{gather*}
B_{n,p}(m_1,\ldots,m_k)=\{(f_1,f_2)|f_i\in C(D^2\rightarrow
M_n(\mathbb{C})),\\
f_1(z)=V_p^*(z)f_2(z)V_p(z),\ |z|=1,\ f_1(0)\in
M_{m_1}(\mathbb{C})\oplus\ldots\oplus M_{m_k}(\mathbb{C})\},
\end{gather*}
where $V_p$ are defined by (\ref{Vk}).

Condition $l_1\equiv l_2\ (mod (m_1,\ldots,m_k)\ )$ implies that
$$\exists c_1,\ldots,c_k\in
\mathbb{Z}:c_1m_1+\ldots+c_km_k=l_1-l_2.$$

Let $H(t)=e^{i 2 \pi c_1 t}E_{m_1}\oplus\ldots\oplus e^{i 2 \pi
c_k t}E_{m_k}\in U(n)$, where $E_m$ denotes the $m\times
m$-identity matrix. Construct the map
$\overline{H}:D^2\longrightarrow U(n)$ by the rule :
\begin{gather*}
z\mapsto H\left(\frac z {|z|}\right),\ z\neq 0,\\
0\mapsto \mathbf{1}_n.
\end{gather*}

Nevertheless $\overline{H}$ is not continuous in $0$, it is easy
to check that $\overline{H}$ defines an isomorphism :
\begin{gather*}
B_{n,l_1}(m_1,\ldots,m_k)\ni(f_1(z),f_2(z))\mapsto
(H^*(z)f_1(z)H(z),f_2)\in \Gamma(\mathcal{B}'_{n,l_2}).
\end{gather*}
The bundle $\mathcal{B}'_{n,l_2}$ is defined by the sewing map
$V'_{l_2}(z)=H^*(z)V_{l_1}(z)$. The equality $ind\ V'_{l_2}=ind\
V_{l_2}$ gives an isomorphism (\ref{bnl12}).
\end{proof}

\begin{rem}
The converse statement to the Theorem 2 takes place. We will not
use it below, so we don't give proof.
\end{rem}

\begin{rem}
One can easily prove the following generalization of theorem 2: if
we "fix" sections of bundles $\mathcal{B}_{n,l_1}$ and
$\mathcal{B}_{n,l_2}$ in a few points in block matrices (with
equal dimension of blocks for both bundles in every point) and
suppose that all blocks are of dimensions $m_1,\ldots,m_k$, then
$l_1\equiv l_2\ (mod (m_1,\ldots,m_k)\ )$ implies isomorphism of
corresponding algebras of sections.
\end{rem}

Let us return to our examples of finitely generated
$F_{2n}$-$C^*$-algebras (see preliminaries).
\begin{prop}
For $C^*$-algebras associated with extended Dynkin diagrams one
has the next isomorphisms:
\begin{align*}
C^*(\widetilde{D_4})&\simeq A_2(1,1;\ 1,1;\ 1,1),\\
C^*(\widetilde{E_6})&\simeq A_3(2,1;\ 1,1,1;\ 1,1,1),\\
C^*(\widetilde{E_7})&\simeq A_4(2,2;\ 2,1,1;\ 1,1,1,1),\\
C^*(\widetilde{E_8})&\simeq A_6(3,3;\ 2,2,1,1;\ 1,1,1,1,1,1).
\end{align*}
\end{prop}

\begin{proof}
Let $A$ denotes $C^*(\widetilde{D_4})$. We can consider this
algebra to be realized by functions on the $S^2$ with values in
matrices of corresponding dimension, not necessary continuous yet.
We can find a covering of $S^2$ by closed subsets $W_j\simeq D^2$
such that every $W_j$ contains only representations in main
dimension or contains one point with reducible representation of
main dimension. Consider the "restrictions" $A(W_j)$ of algebra
$A$ onto subsets $W_j$ (i.e. factor-algebras $A/I(W_j)$, where
$I(W_j)=\set{f\in A,f|_{W_j}\equiv 0}$), then use Theorem 1 to
realize every $A(W_j)$ as sections of some algebraic bundle, i.e.
everyone of $A(W_j)$ coincides with one of the
$B_{n,l}(m_1,\ldots,m_k)$ (possibly with $m_1=n,\
m_2=\ldots=m_k=0$).

Different $A(W_j)$ are "sewed" by some continuous maps, which we
use to construct the bundle $\mathcal{B}_{n,l}$ such that algebra
$A$ is isomorphic to $B_{n,l}(m_1,\ldots,m_k)$. Theorem 2 shows
that we can suppose $l=0$.
\end{proof}

\begin{prop}
Enveloping $C^*$-algebra for $P_{\alpha,\beta}$ is :
$$C^*(P_{\alpha,\beta})=\{f\in C(S^2\longrightarrow
M_2(\mathbb{C}))|f(x_1),\ f(x_2)\in M_1(\mathbb{C})\oplus
M_1(\mathbb{C}),\ x_1,\ x_2\in S^2\}$$
\end{prop}
\begin{proof}
See the proof of previous theorem.
\end{proof}

The following example of enveloping $C^*$-algebra is considered in
\cite{vs}. One can get the result of \cite{vs} using our Theorem
1.
\begin{prop}
Enveloping $C^*$-algebra $\mathcal{P}_2$ for free pair of
projections (see \ref{free}) is :
$$\mathcal{P}_2=\{f\in C([0,1]\longrightarrow
M_2(\mathbb{C}))| f(0),\ f(1)\ \mbox{are diagonal }\}.$$
\end{prop}

\begin{proof}
As in previous case, we realize this algebra by sections of some
bundle over $[0,1]$. By statement 1 all bundles over $[0,1]$ are
trivial.
\end{proof}

All the constructions above in the paper, including Statement
\ref{vas}, can be generalized to prove the next proposition: 

\begin{prop}
Enveloping $C^*$-algebra for $\overline{P}_\eps$ $($example
\ref{com}$)$ has the form:
\begin{gather*}
\overline{P}_e=\{f\in C(S^2\times I\longrightarrow
M_2(\mathbb{C}))|f(x_1,1/2),\ f(x_2,t),\ f(x_3,t)\in
M_1(\mathbb{C})\oplus M_1(\mathbb{C}),\\
t\in I,\ x_1,\ x_2,\ x_3\in S^2\}.
\end{gather*}
\end{prop}

\begin{proof}
Here we use the fact that the theory of algebraic bundles over
$S^2\times I$ is the same as for $S^2$. Then the proof of Theorem
2 should be transferred from $S^2$ to $S^2\times I$.
\end{proof}

\begin{rem}
The definition of triviality for homogeneous $C^*$-algebras can be
generalized. One can define $F_{2N}$-$C^*$-algebra satisfying
conditions of the Statement \ref{vas} to be trivial iff it is
isomorphic to the $F_{2N}$-$C^*$-algebra with the trivial
corresponding bundle
$\mathcal{B}(X,M_N(\mathbb{C}),\underline{U(N)})$. Then all
propositions above state the triviality of $C^*$-algebras under
consideration. Example of non-trivial non-homogeneous
$C^*$-algebra one can construct by "fixing" sections of bundle
$\mathcal{B}_{4,1}$ in the $2\times 2$-block matrices in one
point, it is also an example of algebra that is not defined by the
dual space.
\end{rem}
\noindent
{\bf{Acknowledgements.}}\\
Author expresses his gratitude to Prof. Yuri\u{\i} S.
Samo\u{\i}lenko for the problem stating and permanent attention. I
am also indebted to Dr. Daniil P. Proskurin for help during the
preparation of this paper.

\end{document}